\documentclass[a4paper,12pt,reqno]{amsart}
\usepackage{amssymb}
\usepackage{amsmath}

\def\R{\mathbb{R}}
\def\C{\mathbb{C}}

\def\so{\mathfrak{so}}

\def\SO{\mathbf{SO}}
\def\Spin{\mathbf{Spin}}

\def\G{\mathbf{G}}
\def\diag{\mathrm{diag}\,}
\def\Ric{\mathrm{Ric}}
\def\tr{\mathrm{tr}}
\def\la{\langle}
\def\ra{\rangle}

\def\spin{\mathfrak{spin}}

\def\Cas{\mathrm{Cas}}
\def\End{\mathrm{End}\,}
\def\Hom{\mathrm{Hom}\,}
\def\Hol{\mathrm{Hol}}

\def\id{\mathrm{id}}

\def\pr{\mathrm{pr}}
\def\Aut{\mathrm{Aut}}
\def\ins{\,\lrcorner\,\,}

\def\beq{\begin{equation}}
\def\eeq{\end{equation}}
\def\bea{\begin{eqnarray*}}
\def\eea{\end{eqnarray*}}
\def\ba{\begin{array}}
\def\ea{\end{array}}
\def\proof{\noindent\textbf{Proof:}\quad}
\def\qed{\quad\hfill\ensuremath{\Box}}

\def\L{\Lambda}

\def\r{\end{proof}}

\def\g{\mathfrak{g}}

\def\End{{\rm End}}

\def\d{{d^*}}

\def\es{\,\lrcorner\,}
\def\we{\wedge}

\def\Ric{\mathrm{Ric}}

\def\id{\mathrm{id}}
\def\be{\begin{equation}}

\def\ee{\end{equation}}

\def\pr{\rm{pr}}
\def\tr{\mathrm{tr}}

\def\Hol{\mathrm{Hol}}

\def\so{\mathfrak{so}}
\def\R{\mathbb{R}}

\def\al{i\,\def\diag{\mathrm{diag}\,}}
\def\W{B}
\def\S{\mathrm{Sym}}
\def\#{\sharp}

\newtheorem{ede}{Definition}[section]

\newtheorem{epr}[ede]{Proposition}

\newtheorem{ath}[ede]{Theorem}

\newtheorem{elem}[ede]{Lemma}

\newtheorem{ere}[ede]{Remark}

\newtheorem{ecor}[ede]{Corollary}
 
\address{Uwe Semmelmann\\ Fachbereich Mathematik, Universit{\"a}t
Hamburg\\ Bundesstr. 55,  D-20146 Hamburg, Germany}
\email{Uwe.Semmelmann@math.uni-hamburg.de}

\begin{document}

\title{Killing forms on $\G_2$-- and $\Spin_7$--manifolds}
\author{Uwe Semmelmann}
\thanks{The author would like to thank the Centre de
Math{\'e}matiques de l'Ecole Polytechnique for hospitality during the
preparation of this work.}
\begin{abstract}
Killing forms on Riemannian manifolds are differential forms
whose covariant derivative is totally skew-symmetric. We
prove that on a compact manifold with holonomy $G_2$
or $Spin_7$ any Killing form has to be parallel.
The main tool is a universal Weitzenb{\"o}ck formula.
We show how such a formula can be obtain for any 
given holonomy group and any representation defining 
a vector bundle.
\vspace{0.1cm}

\end{abstract}

\maketitle

\section{Introduction}

Killing forms are a natural generalization of Killing vector fields.
They are defined as differential forms $u$, such that
$\nabla u$ is totally skew-symmetric. More generally one considers
twistor forms,  as forms in the kernel of an elliptic differential 
operator, defined similar to the twistor operator in spin geometry. 
Twistor 1-forms are dual to conformal vector fields.

The notion of Killing forms was introduced by K.~Yano in \cite{yano},
where he already noted that a $p$--form $u$ is a Killing form
if and only if for any geodesic $\gamma$ the $(p-1)$--form 
$\,{\dot \gamma} \,\lrcorner\, u\,$ is parallel along $\gamma$. 
In particular, Killing forms define quadratic first integrals of the
geodesic equation, i.e. functions which are constant along geodesics. 
This motivated an intense study of Killing forms in the physics
literature,~e.g. in the article~\cite{penrose} of R.~Penrose and M.~Walker.
More recently Killing  and twistor forms have been
successfully applied to define symmetries of field equations
(c.f.~\cite{be1}, \cite{be2}).

On the standard sphere the space of twistor forms coincides with
the eigenspace of the Laplace operator for the minimal eigenvalue and
Killing forms are the coclosed minimal eigenforms. The sphere also realizes
the maximal possible number of twistor or Killing forms.
There are only very few
further examples of compact manifolds admitting Killing $p$-forms
with $p\ge 2$. These are Sasakian, nearly K{\"a}hler and weak-$G_2$ manifolds,
and products of them.

The present article is the last step in the study of Killing forms
on manifolds with restricted holonomy. It was already known that
on compact K{\"a}hler manifolds Killing $p$-forms with $p\ge 2$ are
parallel (\cite{yama3}). Moreover we showed in \cite{au1} and
\cite{fau} that the same is true on compact quaternion-K{\"a}hler
manifolds and compact symmetric spaces. Here we will prove the
corresponding statement for the remaining holonomies $\G_2$
and $\Spin_7$.

%%%%%%%%%%%%%%%%%%%%%%%%%%%%%
\begin{ath}\label{main1}
Let $(M^7,\,g)$ be a compact manifold with holonomy  $\G_2$. Then 
any Killing form and any $\ast$--Killing form is parallel.
Moreover, any twistor $p$--form, with $\,p\ne 3,4$, is parallel.
\end{ath}
%%%%%%%%%%%%%%%%%%%%%%%%%%%%%

%%%%%%%%%%%%%%%%%%%%%%%%%%%%%
\begin{ath}\label{main2}
Let $(M^8,\,g)$ be a compact manifold with holonomy  $\Spin_7$. Then 
any Killing form and any $\ast$--Killing form is parallel.
Moreover, any twistor $p$-form, with $\,p\neq 3,4,5$, is parallel.
\end{ath}
%%%%%%%%%%%%%%%%%%%%%%%%%%%%%

The main tool for proving the two theorems are suitable 
Weitzenb{\"o}ck formulas for the irreducible components of
the form bundle. More generally we prove a universal Weitzenb{\"o}ck
formula,~i.e. we show how to obtain for any fixed holonomy
group $G$ and any irreducible $G$-representation $\pi$, 
a Weitzenb{\"o}ck formula
for the twistor operators acting on sections of the vector bundle
defined by $\pi$. Our formula is already known in  the
case of Riemannian holonomy $\SO_n$ (c.f.~\cite{pg1}).
However, it seems to be new and so far unused in the
case of the exceptional holonomies $\G_2$ and $\Spin_7$.
We describe here an approach to Weitzenb{\"o}ck formulas
which is further developed and completed in~\cite{gu}.

{\it Acknowledgments.}
The author would like to thank Andrei Moroianu and Gregor
Weingart for many valuable discussions and the continued
interest in twistor forms. In particular he is grateful
to Gregor Weingart for sharing his views on Weitzenb{\"o}ck
formulas and his helpful comments on the Casimir normalization.

\section{Twistor forms on Riemannian manifolds}\label{eins}

In this section we recall the definition and basic facts on
twistor and Killing forms. More details and further 
references can be found in~\cite{uwe}. Most important for
the later application will be the integrability condition
given in Proposition~\ref{integrabl}.

Consider a $n$--dimensional Euclidean vector space 
$( V,\,\la\cdot,\cdot\ra)$. Then the  tensor product  
$V^*\otimes\Lambda^pV^*$ has the following
$O(n)$--invariant decomposition:
\bea\label{deco3}
V^*\otimes\Lambda^pV^*
&\cong&
\Lambda^{p-1}V^* 
\oplus \Lambda^{p+1}V^* 
\oplus 
\Lambda^{p,1}V^* 
\eea
where $\Lambda^{p,1}V^*$ is the intersection of the kernels of wedge
and inner product. This decomposition  immediately
translates to Riemannian manifolds $\,(M^n,\,g)$, where we have 
\begin{equation}\label{deco1}
T^*M\otimes\Lambda^pT^*M
\;\cong\;
\Lambda^{p-1}T^*M
\oplus
\Lambda^{p+1}T^*M
\oplus
\Lambda^{p,1}T^*M
\end{equation}
with $\Lambda^{p,1}T^*M$ denoting the vector bundle corresponding to 
the representation
 $\Lambda^{p,1}$. The covariant derivative $\nabla \psi$ 
of a $p$--form $\psi$ 
is a section of $\;T^*M\otimes\Lambda^pT^*M$. Its projections onto
the summands $\,\Lambda^{p+1}T^*M\,$ and $\,\Lambda^{p-1}T^*M\,$
are just the differential
$d\psi$ and the codifferential $\d \psi$.
Its projection onto the third summand $\,\Lambda^{p,1}T^*M\,$ defines a
natural first order differential operator $T$, called the 
{\it twistor operator}. The twistor operator
$
T:\Gamma(\Lambda^p T^*M)\,\rightarrow \,\Gamma(\Lambda^{p,1}T^*M) 
\subset 
\Gamma(T^*M\otimes\Lambda^pT^*M)
$
is given for any vector field $X$ by the following formula
\beq\label{twistor}
[\,T\psi\,]\,(X)
\;:=\;
[\pr_{\Lambda^{p,1}}(\nabla \psi)]\,( X)
\;=\;
\nabla_X\, \psi
\;-\;
\tfrac{1}{  p+1}\, X \es d\psi
\;+\;
\tfrac{1}{ n-p+1}\, X^*\,\wedge\,\d \psi \ .
\eeq
From now on we will identify $TM$ with $T^*M$ using the metric.

%%%%%%%%%%%%%%%%
\begin{ede}
A $p$--form $\,\psi\,$ is called a {\it twistor $p$--form} if and only if 
$\,\psi\,$ is in the kernel of $\,T$,~i.e. 
if and only if $\,\psi\,$ satisfies 
\begin{equation}\label{killing}
\nabla_X\,\psi\;=\;
\tfrac{1}{  p+1}\,X\es d\psi \;-\;
\tfrac{1}{  n-p+1}\, X\,\wedge\,\d \psi\ ,
\end{equation}
for all vector fields $X$. 
If the $p$--form $\psi$ is in addition coclosed, it is called a 
{\it Killing $p$--form}. 
%This is equivalent to $\nabla\psi\in\Gamma(\Lambda^{p+1}T^*M)$ or to 
%$X \es \nabla_X \psi = 0$ for any vector field $X$.
A closed twistor form is called {\it $\ast$-Killing form}.
\end{ede}
%%%%%%%%%%%%%%%%

Twistor forms are also known as {\it conformal Killing forms}
or skew-symmetric {\it Killing-Yano tensors}. Twistor 1-forms
are dual to conformal vector fields and Killing 1-forms are
dual to Killing vector fields. Note that the Hodge star-operator 
$\ast$ maps twistor $p$--forms into twistor $(n-p)$--forms. 
In particular, it interchanges Killing and  $\ast$-Killing forms.

Twistor forms are well understood
on compact K{\"a}hler manifolds (c.f.~\cite{au1}). Here they
are closely related to Hamiltonian 2-forms  recently studied in
\cite{pg3}. In particular, one has examples on the complex
projective space in any even degree.

Differentiating Equation~(\ref{twistor}) one obtains the 
two equations
\begin{eqnarray}
\nabla^*\nabla\psi 
&=& \label{weizen1}
\tfrac{1}{ p+1}\,
\d d\, \psi\;\;+\;\; \tfrac{1}{  n-p+1}\,d \d \, \psi\;\;+\;\;T^*T\,
\psi \ ,
\\[1.5ex]
q(R)\,  \psi 
&=&\label{weizen2}
\tfrac{p}{ p+1}\,\d d\, \psi
\;\;+\;\;\tfrac{n-p}{ n-p+1}\,d \d \, \psi\;\;-\;\;T^*T\, \psi \ ,
\end{eqnarray}
where $q(R)$ is the curvature term appearing in the classical
Weitzenb{\"o}ck formula for the Laplacian on $p$--forms:
$\Delta  = \d d\,+\,d\d  = \nabla^*\nabla \,+\, q(R) $. It is  the symmetric
endomorphism of the bundle of differential forms defined by
\begin{equation}\label{qr}
q(R)\;=\;\sum\,e_j\,\wedge\,e_i \,\es \, R_{e_i,e_j},
\end{equation}
where $\,\{e_i\}\,$ is any local orthonormal frame and $R_{e_i,e_j}$
denotes the curvature of the form bundle. On 1-forms the
endomorphism $q(R)$ is just the Ricci curvature.
It is important to note that one may define $q(R)$ also in a more general
context. For this we first rewrite equation~(\ref{qr}) as
$$
q(R)
\;=\;
\sum_{i < j}\,
(e_j\,\wedge \,e_i \es \,-\; e_i\,\wedge \,e_j\es)\, R_{e_i,e_j}
\;=\;
\sum_{i < j}\,
(e_i \we e_j)\, R(e_i \we e_j )\,
$$
where the Riemannian curvature $R$ is considered as element of
$\, \S^2(\L^2 TM)\,$ and 2-forms act via the standard representation
of the Lie algebra $\,\so(T_mM)\cong\L^2 T_mM \,$ on the space of
$p$--forms. 
Note that we can replace $\,\{e_i \we e_j\}\,$ by  any basis  of $\so(T_mM)$
orthonormal with respect to the scalar product induced by $g$ on
$\,\so(T_mM)\cong\L^2 T_mM $.

Let $(M, g)$ be a Riemannian manifold with holonomy group $G=\Hol$. Then
the curvature tensor takes values in the Lie algebra 
$\,\g \, $ of the holonomy group and we can write $q(R)$ as 
$$
q(R)
\;=\;
\sum \, X_i \, R(X_i)\, \qquad \in \;\S^2(\g)
$$
where $\{X_i\}$ is any orthonormal basis of $\, \g \,$ 
acting via the form representation restricted to the holonomy group.
It is clear that in this way $q(R)$ gives rise to a symmetric endomorphism
on any associated vector bundle defined via a representation of
the holonomy group. Moreover this bundle endomorphism
preserves any parallel subbundle and  its action only depends on the
representation defining the subbundle and
not on the particular realization.

Integrating Equation~(\ref{weizen2}) yields
a characterization of twistor forms on compact manifolds. This
generalizes the characterization of Killing vector fields on compact
manifolds, as divergence free vector fields
in the kernel of $\Delta - 2\Ric$.
%%%%%%%%%%%%%%%%%
\begin{epr}\label{integrabl}
Let $(M^n,\,g)$ a compact Riemannian manifold. Then a $p$--form $\psi$ 
is a twistor $p$--form, if and only if
$\,\,
q(R)\,  \psi 
=
\tfrac{p}{ p+1}\,\d d\, \psi\, + \,\tfrac{n-p}{ n-p+1}\,d \d \, \psi \ .
$
A coclosed $p$--form $\psi$ a Killing form if and only if
$\,\,
\nabla^*\nabla \psi = \frac{1}{p}q(R)\psi \ .
$
\end{epr}
%%%%%%%%%%%%%%%%%%%%

One has similar characterizations for closed twistor forms and for
twistor m--forms on 2m--dimensional manifolds.
For the later application in the case of compact Ricci-flat
manifolds we still mention an immediate consequence
of Equation~(\ref{weizen2}).

\begin{ecor}\label{ricci}
Let $M$ be a compact manifold and let $EM \subset \Lambda^pT^*M$ 
be a parallel subbundle such that $\,q(R)\,$ acts trivially on $\,E$. 
Then any twistor and any harmonic form in
$\,EM$ has to be parallel. In particular, there are only parallel
twistor resp. harmonic forms in $\, \Lambda^k_7\,$ for \,$1\le k \le
6\,$  on a 
compact $\,\G_2$-manifolds and there are only parallel 
twistor resp.  harmonic forms in
$\,\Lambda^k_7\,$ for $\,k=2,4,6\,$ and in  $\,\Lambda^k_8\,$ for 
$\,k=1,3,4,5,7\,$ on a compact $\,\Spin_7$-manifold.
\end{ecor}
\proof
It remains to prove that $q(R)$ acts trivially on the
given bundles. This is clear for the 7-dimensional bundles in
the $\G_2$ case and for the 8-dimensional bundles in
the $\Spin_7$-case. Indeed, these bundles are defined
via the holonomy representation, where $q(R)$ acts as 
Ricci-curvature, which is zero in the case of
holonomy  $\G_2$ or  $\Spin_7$.

In the $\Spin_7$-case there is still an other argument
for proving the trivial action of $q(R)$ on the
8- and 7-dimensional part:
The spinor bundle of a manifold with $\Spin_7$--holonomy
splits into the sum of a trivial line bundle, corresponding to the
parallel spinor, and the sum of a bundle of rank 7 and a bundle of rank 8. 
These two bundles are induced by the 8-dimensional holonomy
representation and by the 7-dimensional standard representation. 
It is well-known that $\,q(R)\,$ acts as
$\,\tfrac{s}{16}\id\,$ on the summands of the spinor bundle. But for
$\,\Spin_7$--manifolds the scalar curvature $\,s\,$ is zero and
we conclude that $\,q(R)\,$ acts trivially on the bundles
in question.
\qed

\section{A universal Weitzenb{\"o}ck formula}

In this section we derive for any manifold with a 
fixed holonomy  one basic Weitzenb{\"o}ck
formula. The coefficients of this formula will
depend on the holonomy group and the defining
representation. This is only the first step of a more 
general method of producing all possible
Weitzenb{\"o}ck formulas on such manifolds
(c.f.~\cite{gu}).  

We consider the following situation: let  $(M^n, g)$ be 
an oriented Riemannian manifold with holonomy group 
$G:= \Hol(M,g) \subset \SO(n)$. Then the $\SO(n)$-frame bundle
reduces to a $G$-principal bundle $P_G\rightarrow M$
and all natural bundles over $M$ are associated to
$P_G$ via representations of $G$. 

If $\pi: G \rightarrow \Aut(E)$ is a complex representation of $G$
we denote with $EM$ the corresponding associated bundle
over $M$. In particular, we will denote the complexified holonomy 
representation of $G$, given by the inclusion $G \subset \SO(n)$,
with $T$. The associated vector bundle is then of course
the complexified tangent bundle $TM$.

The Levi-Civita connection of $(M,g)$ induces a connection $\nabla$ on
any bundle $EM$ and the covariant derivative of a section of $EM$ is
a section of $TM \otimes EM$. Hence, we may define natural first order 
differential operators by composing the covariant derivative with
projections onto the components of $TM \otimes EM$. These
operators are also known as {\it Stein-Weiss operators}.

Let $T\otimes E = \sum E_\al$ be the decomposition of $T\otimes E$
into irreducible $G$-representations, where we consider 
$E_\al$ as a subspace of $T\otimes E$. This induces a
corresponding decomposition of the tensor product 
$TM \otimes EM$. For any component $E_\al M$ we define a 
{\it twistor operator} $T_\al$ by
$$
T_\al: \Gamma(EM) \rightarrow \Gamma(E_\al M), 
\qquad
T_\al (\psi):= \pr_\al (\nabla \psi) \ ,
$$
where $\pr_\al$ denotes the projection $T \otimes E \rightarrow
E_\al  \subset T \otimes E $ and the corresponding bundle map. 
In the  following
we will make no difference between representations 
resp. equivariant maps and the corresponding vector bundles resp.
bundle maps.

Since we are on a Riemannian manifold we have for
any twistor operator $T_\al$ its formal adjoint 
$T^*_\al: \Gamma(E_\al M)\rightarrow \Gamma(EM)$.
The aim of this section is to derive a  Weitzenb{\"o}ck
formula for the second order operators $\,T^*_\al\circ T_\al$,~i.e.
a linear combination
$\;
\sum_\al \,c_\al \, T^*_\al\circ T_\al
$\;
with real numbers $c_\al $, which is of zero order,~i.e.
a curvature term. The coefficients  $c_\al $ will depend
on the holonomy group $G$ and the representation $E$.

Our approach to Weitzenb{\"o}ck formulas, further developed in 
\cite{gu}, is motivated by the
following remarks. Let $\psi$ be any section of $EM$, then
$\nabla^2 \psi$ is a section of the bundle $TM \otimes TM \otimes EM$.
Any $G$-equivariant homomorphism 
$F \in \Hom_G(T \otimes T \otimes E, E)$ defines by
$\psi \mapsto F(\nabla^2\psi)$  a second order differential
operator acting on sections of $EM$. For describing these 
homomorphisms it is rather helpful to use the natural
identification
\medskip
$$
\Hom_G(T \otimes T \otimes E, E)
\;\cong\;
\End_G(T \otimes E)
\;\cong\;
\Hom_G(T \otimes T, \End E) \ .
$$

A homomorphisms $F: T\otimes T \rightarrow \End E$
is mapped onto the endomorphism $F$ of $T \otimes E$
defined by
$\,
F(a \otimes s) = \sum e_i \otimes F_{e_i \otimes a}(s)
$,
for any orthonormal basis $\{e_i\}$ of $T$ and any $a\in T, s\in E$. 
Conversely an endomorphism $F$ is mapped to the homomorphism $F$ 
with $F_{a\otimes b} (s)= a \ins F(b\otimes s)$. 
In particular, the identity of $T\otimes E$ 
is mapped onto $\id_{a \otimes b} =  g(a,\,b)\,\id_E$.
Finally, $F\in \Hom(T \otimes T \otimes E, E) $ is
defined as $F(a\otimes b\otimes s)= F_{a \otimes b} s$.

Beside the identity $\,\id_{T \otimes E}\,$ we have the projections 
$\,\pr_\al: 
T \otimes E \rightarrow E_\al \subset T \otimes E
$\,
as important examples of invariant endomorphisms.
The following proposition describes the corresponding second order 
differential operators.

\begin{epr}
Let $T\otimes E = \sum E_\al$ the decomposition into
irreducible summands, with corresponding twistor operators
$T_\al$. Then any section $\,\psi\,$ of $\,EM\,$ satisfies
$$
\begin{array}{ll}
&(1)
\qquad \id (\nabla^2 \psi) 
\,\quad=\,\, 
- \, \nabla^*\nabla \psi
\\[1.ex]
&(2)
\qquad \pr_\al (\nabla^2 \psi) 
\,\,\,=\,\, 
-\,T^*_\al \,\,T^{\phantom{*}}_\al (\psi)
\end{array}
$$
\end{epr}
\proof
Let $\{e_i\}$ be a parallel local ortho-normal frame of $TM$, then
$\nabla^2 = \sum e_i \otimes e_j \otimes \nabla_{e_i}\nabla_{e_j}$
and
$\,
\id (\nabla^2\psi)
\,=\,
\sum g(e_i, e_j) \nabla_{e_i}\nabla_{e_j}\psi
\,=\,
-\, \nabla^*\nabla \psi \ ,
$
which proves Equation~(1). 

Equation~(2)  is a direct
consequence of the following more general statement.

%%%%%%%%%%%%%%%%%%%
\begin{elem}\label{adjoint}
Let $(M,\,g)$ be a Riemannian manifold and let $E,\,F$
be hermitian vector bundles over $M$, equipped with
metric connections $\nabla$. If 
$D:\Gamma(E)\rightarrow \Gamma(F)$ is a differential
operator defined as $\,D = p \circ \nabla$, where 
$\,p : TM \otimes E \rightarrow F\,$ is some parallel
linear map. Then the adjoint operator for $D$ is
$\,D^* = \nabla^* \circ p^*\,$ and
$\,
D^*D 
\;=\;
-\,\tr\circ (\id \otimes p^*p)\,\nabla^2 \ .
$
\end{elem}
%%%%%%%%%%%%%%%%%%%%%%

The second equation follows with
$F= TM \otimes EM$  and an orthogonal projection $\,p=\pr_i$
 onto a subbundle $\, EM_i \subset TM \otimes EM$.
Indeed, in this case we have $p^*p = p^2 = p$ and it follows
$\,(T_i) T_i\psi = - \tr \circ (\id \otimes \pr_i) \nabla^2\psi
= - \pr_i\,(\nabla^2\,\psi)
$.

We will now prove the lemma. Note, that the formal 
adjoint $\,\nabla^*\,$ of the covariant derivative
$\nabla:\Gamma(E)\rightarrow
\Gamma( TM \otimes E )\,$ is given as the composition
of the following differential operators:
$\,
 \Gamma( TM \otimes E )
\stackrel{\nabla}
   \longrightarrow
 \Gamma( TM \otimes TM\otimes E )
 \stackrel{-\textrm{tr}\,\,}
   \longrightarrow
 \Gamma( E ) \ ,
$
where $\,\nabla\,$ also denotes the tensor product connection, i.e 
$\,
\nabla
:=
\nabla\otimes\textrm{id}\,+\,\textrm{id}\otimes\nabla
$.
Indeed, we obtain for a vector field  $X\in\Gamma( TM )$ 
and a section $\psi\in\Gamma( E )$ that
$$
 X\otimes\psi
 \longmapsto
 \sum_i
 \left( e_i \otimes \nabla_{e_i}X \otimes \psi
        + e_i \otimes X \otimes \nabla_{e_i} \psi \right)
 \longmapsto
 \sum_i
(- \mathrm{div}(X) \,\psi - \nabla_{X}\psi )
$$
and it is easily seen that this composition is formally adjoint
to $\,\nabla$. 
Since $D$ is defined as $D=p \circ \nabla$ we have
$D^* = (p \circ \nabla)^*=\nabla^* \circ p^*$,
thus
$\,
D^*\,D
=
(p \circ \nabla)^*  (p \circ \nabla)
=
\nabla^* \circ p^*p \circ \nabla \ .
$\,
Since $\,p\,$ is a parallel map it commutes with
$\nabla^*$. Hence we can substitute the formula
for $\nabla^*$ to obtain
$
D^*D
=
\nabla^* \circ p^*p \circ \nabla
=
 - \,\textrm{tr} \circ \nabla \circ p^*p \circ \nabla
=
 -\, \textrm{tr} \circ (\textrm{id} \otimes p^*p) \circ \nabla^2 \ . 
$

\qed

Since we obviously have $\id = \sum \pr_\al $, the proposition
above immediately implies  a rather useful
formula for the operator $ \nabla^*\nabla $, which 
corresponds to  Equation~(\ref{weizen1}) in the
case $G=\SO_n$ and $E= \Lambda^pT$.

\medskip

\begin{ecor}
$\qquad\qquad\qquad\qquad
\nabla^*\nabla \;\;=\;\; \sum_\al  \, T^*_\al\circ T_\al
$
\end{ecor}

Let $G$ be the holonomy group of an irreducible, non-symmetric 
Riemannian manifold. It is then well-known that any
isotypic component of $T\otimes E$ is irreducible,~i.e.
in the decomposition $T\otimes E = \sum E_\al$ any
summand $E_\al$ occurs only once. As a consequence 
 the projection maps $\, \{ \pr_\al \} \,$ form
a basis of  $\End_G(T \otimes E)$ and any invariant
endomorphism $F$ of $T \otimes E$ can be written as
$F= \sum f_\al \pr_\al$, with $F|_{E_\al} = f_\al \id$.

It turns out that for certain invariant endomorphisms
$F$ the operator $F\circ \nabla^2$ is in fact a zero order
term. Hence, in these cases $F$ gives rise to the
Weitzenb{\"o}ck formula 
$\, F\circ \nabla^2\,=\,\sum\, f_\al \,T^*_\al\, T_\al  $.
The following lemma will provides  us with an easy criterion for
deciding which invariant  endomorphisms $F$ lead
to Weitzenb{\"o}ck formulas.

\begin{elem}
Let $F$ be an equivariant endomorphism of $\,T\otimes E\,$
considered as element of $\,\Hom_G(T \otimes T, \End E)$.
Then $\,F \circ \nabla^2\,$ defines a zero order operator if and only
if $ \,F_{a \otimes b} = - F_{b\otimes a}\,$ for any 
vectors $\,a, b \in T$.
\end{elem}
\proof
Let $R$ be the curvature of $EM$ and
let $\{e_i\}$ be a parallel local frame, then
$$
F \circ \nabla^2
=
\sum F(e_i \otimes e_j)\nabla_{e_i} \nabla_{e_j}
=
\tfrac{1}{2} \sum F(e_i \otimes e_j)
(\nabla_{e_i} \nabla_{e_j} - \nabla_{e_j} \nabla_{e_i})
=
\tfrac{1}{2} \sum F(e_i \otimes e_j) R_{e_i, e_j}
$$
\qed

We show in \cite{gu} that $\End_G(T \otimes E)$ is in many cases,
including the exceptional holo\-nomies $\G_2$ and $\Spin_7$, the
quotient of a
polynomial algebra generated by one special endomorphism,
the {\it conformal weight operator} $\W$. 

\begin{ede}
The conformal weight operator 
$\W \in \End_G(T \otimes E) \cong \Hom_G(T \otimes T, \End E)$
is for any $a, b \in T, s \in E$ defined as
$$
\W_{a \otimes b} \, s 
\,\,:=\,\,
\pr_\g (a \wedge b) \, s \ ,
$$
where $\g$ is the Lie algebra of the holonomy group $G$ and
$\pr_\g$ denotes the projection 
$\Lambda^2T \rightarrow \g \subset \so_n \cong \Lambda^2T $.
Here $\g$ acts via the differential of the representation
$\pi$ on $E$.
\end{ede}

However, for the present article it is only important to note
that $B$ defines a Weitzen\-b{\"o}ck formula, since obviously
$B_{a \otimes b}= - B_{b\otimes a}$, for any $a, b \in T$. We will later apply this
formula for proving that any Killing form on a compact
manifold with exceptional holonomy has to be parallel. 

The curvature term defined by $B$ turns out to be the
endomorphism $q(R)$ already introduced in Section~\ref{eins}.
In fact, Equation~(\ref{weizen2}) can be
considered as the Weitzenb{\"o}ck formula corresponding
to $B$ in the special case of $\,G=\SO_n\,$ and $\,E=\Lambda^pT$.

\begin{elem}
$
\qquad\qquad\qquad\qquad
\W \circ \nabla^2 \quad = \quad q(R)
$
\end{elem}
\proof
Let $\{X_i\}$ be an ortho-normal basis for the induced scalar 
product on $\,\g \subset \Lambda^2T$
and let $\{e_i\}$ be a local ortho-normal frame. Then
\bea
\W\circ \nabla^2
&=&
\sum \,\pr_\g(e_i \wedge e_j)\, \nabla^2_{e_i, e_j}
\,\,=\,\,
\tfrac12\,\sum \pr_\g(e_i \wedge e_j)
( \nabla^2_{e_i, e_j} \,-\,  \nabla^2_{e_j, e_i})
\\[1ex]
&=&
\sum_{i < j}\, \pr_\g(e_i \wedge e_j)\, R_{e_i, e_j}
\,\,=\,\,
\sum\, X_i \cdot R(X_i)\cdot \,\,=\,\,q(R) \ .
\eea
\qed

In order to obtain the general Weitzenb{\"o}ck formula defined
by $B$ we have to write $B = \sum b_\al \pr_\al$ and to
determine the coefficients $b_\al$.  We first describe 
the conformal weight operator as an element of 
$\End_G(T \otimes E)$.

\begin{elem}\label{opB}
Let $\{X_i\}$ be an orthonormal basis for the induced scalar product on 
$\,\g \subset \Lambda^2T$. Then $\,B = - \sum \,X_i \otimes X_i$,
where $X_i$ is acting on $T$ resp. $E$ via the holonomy 
representation resp. the representation $E$.
\end{elem}
\proof
Using the formula $\la X, a\wedge b \ra = \la Xa, b\ra $,
for $a, b \in T $ and $X \in \Lambda^2T \cong \so_n$, 
we may write $B$ as
\bea
\W(a \otimes s) 
&=&
\quad \sum\,e_i \otimes \pr_\g(e_i \wedge a)\, s
\; = \;
\sum\,e_i \otimes \la e_i \wedge a,\,X_j \ra \, X_j \, s
\\[1ex]
&=&
\quad \sum\, \la X_j e_i,\,a\ra\, e_i \otimes  X_j \, s
\; = \;
-\,\sum\,\la  e_i,\,X_ja\ra\, e_i \otimes  X_j \, s
\\[1ex]
&=&
-\,(\sum X_j \otimes X_j)\; a\otimes s
\eea
\qed

Let $G$ be a compact semi--simple Lie group, with Lie algebra 
$\mathfrak g$ and let $\pi:G \rightarrow \Aut(V)$ be a representation 
of $G$ on the complex vector space $V$. 
If $\{X_i\}$ is a  basis of $\g$, orthonormal with respect to an
invariant scalar product $g$, the  
{\it Casimir operator} $\Cas^g_\pi\in \End(V)$ is defined as
$$
\Cas^g_\pi \;:=\; \sum \pi_*(X_i)\circ \pi_*(X_i) \;=\;\sum \, X_i^2
\ ,
$$
where $\pi_*:\g \rightarrow \End(V)$ denotes the differential of the
representation $\pi$. It is well-known that $\Cas^g_\pi = c^g_\pi\, \id_V$,
if the representation $\pi$ is irreducible. Moreover, the Casimir
eigenvalues $c^g_\pi$ can be expressed in terms of the highest weight 
of $\pi$.

It follows from the lemma above that the conformal weight operator 
$\W$ can be written as a linear combination of Casimir operators,
which leads to

%Hence, the coefficients $b_i$ can be expressed in terms of
%Casimir eigenvalues.

\begin{ecor}\label{confW}
Let $\,T \otimes E = \oplus \,E_\al \,$ be the decomposition of
the tensor product into irreducible components. 
Then the conformal weight operator $\,\W\,$ is given as
\beq\label{bi}
\W \,\,\,=\,\,\,
\sum\, 
b_i\,\pr_\al
\qquad\quad\mbox{with}\qquad\quad
b_i\;=\;  \frac{1}{2}\,(c^{\Lambda^2}_T \,+\, c^{\Lambda^2}_E  \,-\, 
c^{\Lambda^2}_{E_\al}) \ ,
\eeq
for Casimir eigenvalues $\,c^{\Lambda^2}_\pi$  computed with
respect to the induced scalar product on $\, \g \subset \Lambda^2T$.
The corresponding universal Weitzenb{\"o}ck formula on sections of
$EM$ is
\beq\label{weizen3}
q(R)\,\,\,=\,\,\, -\,\sum\, b_i\,\, T^*_i \,T_i \ .
\eeq
\end{ecor}
\proof
Expanding the Casimir operator 
$\,\Cas^{\Lambda^2}_{T \otimes E} = \sum X_i^2\,$
acting on $\,T \otimes E \,$ we obtain
$$
\Cas^{\Lambda^2}_{T \otimes E} \,=\,
\sum \, 
( X_i^2 \otimes \id_E \,+\,
2X_i \otimes X_i \,+\,  \id_T \otimes X_i^2 )
$$
Hence, Lemma~\ref{opB} implies that the conformal weight operator
can be written as
$$
\W
\,\,\,=\,\,\,
-\,\frac12\,
(\Cas^{\Lambda^2}_{T \otimes E}\, \,-\,\, \Cas^{\Lambda^2}_T \otimes \id_E 
\,\,-\,\, \id_T \otimes \Cas^{\Lambda^2}_E) \ .
$$
which yields the formula above after restriction to the irreducible
components $E_\al$.
\qed

\begin{ere}
In the case of Riemannian holonomy $G=\SO_n$
the Weitzenb{\"o}ck formula~(\ref{weizen3}) was
considered for the first time in \cite{pg1}. In this article one
can also find the conformal weight operator and its
expression in terms of Casimir operators. The operator
$B$ appears also in~\cite{pg2}.
Similar results can  be found in~\cite{homma}.
\end{ere}

It remains to compute the Casimir eigenvalues. For doing so we first
recall how to compute them for an irreducible representation of
highest weight $\lambda$ and with respect to the scalar product 
$(\cdot, \cdot)$ defined by the Killing form~$B$.
Let $\rho$ be the half sum of the positive roots of $\g$,  then
\beq\label{casimir}
c^B_\pi \;=\; 
\|\rho\|^2 \;-\; \|\lambda + \rho\|^2
 \;=\;
-\, (\lambda,\,\lambda \, + \,2\rho) \ .
\eeq

For the application of Corollary~\ref{confW} we need the 
Casimir eigenvalues $\,c^{\Lambda^2}_\pi$ defined with respect 
to the induced scalar product on $\, \g \subset \Lambda^2T$. 
The relation between these Casimir eigenvalues is contained
in the following normalization lemma.

\begin{elem}\label{casimir2}
Let $\, \g\,$ be the Lie algebra of a compact simple Lie group 
and let $\,V\,$ be an irreducible real $\g$-representation with
invariant scalar product $\,\la \cdot,\cdot\ra$. If $\pi$ is
any other irreducible $\g$-representation with
invariant scalar product $g$, then
$$
c^{\Lambda^2}_\pi
\;=\;
-\,2\,\frac{\dim \g}{\dim V} \; \frac{1}{c^g_V} \;\,  c^g_\pi \ .
$$
%where $c^{\Lambda^2}_\pi$ denotes the Casimir eigenvalue with
%respect to the scalar product induced by  $\,\la \cdot,\cdot\ra$
%on $\,\g \subset \so(T) \cong \Lambda^2T$.
In particular, the  Casimir eigenvalue of the representation 
$V$ is given as
$\,
c^{\Lambda^2}_V = -2 \frac{\dim \g}{\dim V} \ .
$
\end{elem}
\proof
Since we assume $\,V\,$ to be equipped with a $\,\g$--invariant scalar product 
$\,\la \cdot,\cdot\ra$ we have
$\,\g \subset \so(V) \cong \Lambda^2V$. Restricting the 
induced scalar product onto  $\,\g \subset \Lambda^2V\,$ defines
the natural scalar product 
$\,\la \cdot,\cdot\ra_{\Lambda^2}\,$ on $\,\g$.
Note that
$\,
\la \alpha,\, \beta \ra_{\Lambda^2}
=
-\frac12\,\tr_V(\alpha \circ \beta)
=
\frac12\,\la\alpha,\,\beta\ra_{\End V} \ .
$
Let $\{X_a\}$ be an orthonormal basis of $\g$ with respect to
$\la \cdot,\cdot\ra_{\Lambda^2}$ and let $\{e_i\}$ be an
orthonormal basis of $V$. Then
$\,\,
\Cas^{\Lambda^2}_V (v) = c^{\Lambda^2}_V v = \sum_a X^2_a(v)\,\,
$,
for any $v \in V$,
and we obtain
$$
\dim V \, c^{\Lambda^2}_V 
=
\sum_{a,j}\,\la X_a^2(e_j), \,e_j  \ra
=
-\,\sum_{a,j}\,\la X_a(e_j),\, X_a(e_j)  \ra
=
-\,2\,\sum\,|X_a|^2_{\Lambda^2}
=
-\,2\,\dim \g
$$
which proves the lemma in the case $\pi = V$. Since $\g$ is
a simple Lie algebra it follows that two Casimir operators
defined with respect to different scalar products differ only
by a factor independent from the irreducible representation $\pi$. 
Hence
$\,
\frac{c^{\Lambda^2}_\pi}{c^g_\pi}
\,\,=\,\,
\frac{c^{\Lambda^2}_V}{c^g_V}
$
and the statement of the lemma follows from the special case $\pi=V$.
\qed

\bigskip

In the remaining part of this section we will consider the holonomy
groups $\G_2$ and $\Spin_7$ and make the Weitzenb{\"o}ck 
formula~(\ref{weizen3}) explicit for certain representations
appearing in the decomposition of the form spaces.

\subsection{The group $\G_2$}

The group $\G_2 \subset \SO(7)$ is a compact simple Lie group
of dimension 14 and of rank 2. As fundamental weights one
usually considers $\omega_1$ corresponding to the 7-dimensional
holonomy representation $T$ and $\omega_2$ corresponding
to the 14-dimensional adjoint representation $\g_2$. 
The half-sum of positive roots is the sum of the fundamental
weights,~i.e. $\rho = \omega_1 + \omega_2$.
Any other irreducible $\G_2$-representation can be parameterized
as $\Gamma_{a,b} = a \omega_1 + b \omega_2$, e.g. the
trivial representation is $\Gamma_{0,0} = \C$. Further examples are
$$
\Gamma_{0,1} = \Lambda^2_{14} = \g_2,\qquad
\Gamma_{2,0} = \Lambda^3_{27},\qquad
\Gamma_{1,1} = V_{64} ,\qquad
\Gamma_{3,0} = V^-_{77} \ ,
$$
where the subscripts denote the dimension of the representation,
which is unique up to dimension 77. In dimension 77 on has two
irreducibel $\G_2$--representations, denoted by $V^+_{77}$ and $V^-_{77}$.
Below we need the following decomposition of the spaces of
2- and 3-forms
\beq\label{deco4}
\Lambda^2T \;\cong\; \Lambda^5T \;\cong\;T \oplus \Lambda^2_{14},\qquad\qquad
\Lambda^3T \;\cong\; \Lambda^4T \;\cong\; \C \oplus T \oplus \Lambda^3_{27} \ .
\eeq
Since we want to apply the Weitzenb{\"o}ck formula for the
bundles $\Lambda^2_{14}T$ and $\Lambda^3_{27}T$ we still 
need the following tensor product decompositions
\beq\label{deco}
T \otimes \Lambda^2_{14}
\;\cong\;
T \oplus \Lambda^3_{27} \oplus V_{64},\qquad
T \otimes \Lambda^3_{27} 
\;\cong\;
T \oplus \Lambda^4_{27} \oplus \Lambda^2_{14} 
\oplus V_{64} \oplus  V^-_{77}
\eeq
There is a suitable invariant bilinear form $g$ on $\g_2$, which 
induces the scalar products
$$
g(\omega_1, \, \omega_1) \,\,=\,\, 1 ,\qquad 
g(\omega_2, \, \omega_2) \,\,=\,\, 3 ,\qquad
g(\omega_1, \, \omega_2) \,\,=\,\, \tfrac32\ .
$$
Using Equation~(\ref{casimir}) and Lemma \ref{casimir2}
we obtain the following Casimir eigenvalues
$$
c^{\Lambda^2}_{\Gamma_{a,b}}
\,\,=\,\,
-\,\tfrac23\,
c^{g}_{\Gamma_{a,b}}
\,\,=\,\,
-\,\tfrac23\,
(a^2 \,\,+\,\, 3\,b^2 \,\,+\,\, 3\,a\,b \,\,+\,\, 5\,a \,\,+\,\,9\,b)
\ .
$$
In particular we have
$\,\,
c^{\Lambda^2}_{ \Lambda^3_{27}} =  -\,\frac{28}{3},\,\,
c^{\Lambda^2}_{\Lambda^2_{14} } = -\,8,\,\,
c^{\Lambda^2}_{T} = -\,4,\,\,
c^{\Lambda^2}_{V_{64}} = -\,14,\,\,
c^{\Lambda^2}_{V^-_{77}} = -\,16
\ .
$

\noindent
Finally, we use (\ref{bi}) to obtain the Weitzenb{\"o}ck formula
on the bundles  $\,\Lambda^2_{14}\,$ and $\,\Lambda^3_{27}$.
Recall that the twistor operator $T_i$ is the projection of the 
covariant derivative onto the i.th summand in the tensor product 
decomposition of $T\otimes E$,~i.e we will number the operators 
$T_i$ according to the numbering of the summands in this decomposition,
which has to be fixed in order to make the notation unique.

Here we will consider the tensor product decomposition given in~(\ref{deco}),
e.g. in the case of the representation $\Lambda^2_{14}$
the operator $T_3$ denotes the projection of the covariant
derivative onto the  summand $V_{64}$,
whereas for the representation  $\Lambda^3_{27}$ the
operator $T_3$ denotes the projection of the covariant
derivative onto the summand $\Lambda^2_{14}$. 

\begin{epr}\label{final1}
Let $\{T_i\}$ be the twistor operators defined corresponding
to the decompositions in~(\ref{deco}). Then the following
Weitzenb{\"o}ck formulas hold
\bea
on \quad \Lambda^2_{14} :&&
q(R) \,\,=\,\, 4\, T^*_1 T_1 \;+\; \tfrac{4}{3}\,  T^*_2 T_2 \;-\; T^*_3T_3
\\[1ex]
on \quad \Lambda^3_{27}:&&
q(R)  \,\,=\,\, \tfrac{14}{3}\, T^*_1 T_1 \;+\; 2\,  T^*_2 T_2 \;+\; 
\tfrac{8}{3}\,T^*_3T_3\;-\; \tfrac{1}{3}\,T^*_4T_4\;-\;\tfrac{4}{3}\,
T^*_5T_5
\eea 
\end{epr}

\subsection{The group $\Spin_7$}
Let $e_1,e_2,e_3$ be the weights of the 7-dimensional
standard representation of $\Spin_7$. Then the
fundamental weights are defined as
$$
\omega_1 = e_1,\quad
\omega_2 = e_1 + e_2,\quad
\omega_3 = \tfrac{1}{2}(e_1 + e_2 + e_3)
$$
corresponding to the representations 
$\Lambda^1\R^7,\Lambda^2\R^7$ and the spin  representation.
All other irreducible $\Spin_7$-representations are
parameterized as $\Gamma_{a,b,c}= a\omega_1 + b\omega_2 + c\omega_3$.
The half-sum of positive roots is given as
$\,\rho = \omega_1+\omega_2+\omega_3 =\frac{1}{2}(5,3,1)$.

The holonomy group $\Spin_7$ is considered as subgroup
of $\SO_8$ such that the holonomy representation $T$ 
is given by the 8-dimensional real spin representation.
This leads to the following decompositions of the
form spaces $\Lambda^kT$.
\beq\label{deco5}
\Lambda^2T \cong \Lambda^2_7 \oplus \Lambda^2_{21},\qquad
\Lambda^3T \cong \Lambda^3_8\oplus \Lambda^3_{48},\qquad
\Lambda^4T \cong \Lambda^4_{1} \oplus \Lambda^4_{7}
\oplus \Lambda^4_{27} \oplus \Lambda^4_{35}
\eeq
Again, the subscripts denote the dimensions of the representations
and of course we have $\Lambda^2_7 \cong \Lambda^4_7$
and $\Lambda^3_8 =T$. 
For the investigation of twistor forms on $\Spin_7$-manifolds
we need Weitzenb{\"o}ck formulas on the bundles corresponding 
to $\;\Lambda^2_{21},\; \Lambda^3_{48},\;  \Lambda^4_{27}\; $
and $\;  \Lambda^4_{35} $. The decompositions of the 
tensor products $T \otimes E$ are given as
\begin{eqnarray}\label{deco6}\nonumber
T \otimes \Lambda^2_{21} &\cong&
T\;\oplus\; \Lambda^3_{48}\;\oplus\; V^a_{112},\qquad
T \otimes  \Lambda^4_{27} \;\cong\;\Lambda^3_{48}\; \oplus\;
V_{168}
\\[1ex]
T \otimes \Lambda^3_{48} &\cong&
\Lambda^4_{35}\; \oplus\;  \Lambda^2_{21}\; \oplus\; \Lambda^2_{7}
\;\oplus\; \Lambda^4_{27}\; \oplus\; V_{105} \;\oplus\; V_{189}
\\[1ex]
T \otimes \Lambda^4_{35} &\cong&
T\oplus\; \Lambda^3_{48}\;
\oplus\; V^a_{112}\;\oplus\; V^b_{112}
\nonumber
\end{eqnarray}
There are two 112-dimensional irreducible $\Spin_7$-representation,
which we denote with $V^a_{112}$ and $V^b_{112}$. In terms of
fundamental weights the representations appearing in the above
decompositions are given as follows
\bea
&&
\Lambda^2_7 \;=\; \Gamma_{1,0,0},\quad
\Lambda^2_{21} \;=\; \Gamma_{0,1,0},\quad
\Lambda^3_{48}\;=\;\Gamma_{1,0,1},\quad
\Lambda^4_{27} \;=\; \Gamma_{2,0,0},\quad
\Lambda^4_{35} \;=\; \Gamma_{0,0,2}
\\[1ex]
&&
V^a_{112} \;=\; \Gamma_{0,1,1},\quad
V^b_{112} \;=\; \Gamma_{0,0,3},\quad
V_{168} \;=\;\Gamma_{2,0,1},\quad
V_{105} \;=\;\Gamma_{1,1,0},\quad
V_{189} \;=\; \Gamma_{1,0,2}
\eea

Next we have to calculate all the necessary Casimir eigenvalues.
We choose on $\spin_7$ an invariant scalar product $g_0$
inducing the Euclidean scalar product on the Lie algebra of
the maximal torus, which is identified with $\R^3$. 
Let us start with the representation $V=\Gamma_{1,0,0} =\Lambda^2_7$
with highest weight $e_1$. Formula~(\ref{casimir}) yields
the Casimir eigenvalue $c^{g_0}_V=-6$ and from Lemma~\ref{casimir2},
for $T=V$,  we also get $c^{\Lambda^2}_V=-6$. Hence, we
conclude $c^{\Lambda^2}_\pi =c^{g_0}_\pi$ for any 
irreducible $\Spin_7$-representation $\pi$. Using the
standard scalar product on $\R^3$ it is now easy to 
calculate the following Casimir eigenvalues
\bea
&&
c^{\Lambda^2}_{\Lambda^2_{21}} \,=\, -10,\quad
c^{\Lambda^2}_{T} \,=\, -\tfrac{21}{4} ,\quad
c^{\Lambda^2}_{\Lambda^4_{27}}\,=\, -14 ,\quad
c^{\Lambda^2}_{\Lambda^4_{35}}\,=\,  -12,\quad
c^{\Lambda^2}_{\Lambda^3_{48}}\,=\, -\tfrac{49}{4} 
\\[1ex]
&&
c^{\Lambda^2}_{V^a_{112}} \,=\, -\tfrac{69}{4} ,\quad
c^{\Lambda^2}_{V_{105}}\,=\, -18 ,\quad
c^{\Lambda^2}_{V_{189}}\,=\, -20 ,\quad
c^{\Lambda^2}_{V_{168}}\,=\, -\tfrac{85}{4} ,\quad
c^{\Lambda^2}_{V^b_{112}} \,=\, -\tfrac{81}{4} ,\quad
\eea

Finally, we use~(\ref{bi}) to compute the coefficients
of the Weitzenb{\"o}ck formula on the bundles
$\;\Lambda^2_{21},\; \Lambda^3_{48},\;  \Lambda^4_{27}\; $
and $\;  \Lambda^4_{35} $. As in the $\G_2$-case we
number the twistor operators corresponding to the decomposition~(\ref{deco6}),~e.g.
for the representation $\Lambda^2_{21}$ the twistor operator
$T_1$ denotes the projection of the covariant derivative onto the
summand $T$  and for the representation
$\Lambda^3_{48}$ it denotes the projection 
onto the  summand $\Lambda^4_{35}$.

\begin{epr}\label{final2}
Let $\{T_i\}$ be the twistor operators defined corresponding
to the decompositions in~(\ref{deco6}). Then the following
Weitzenb{\"o}ck formulas hold
\bea
on \quad \Lambda^2_{21}:&&
q(R) \,\,=\,\,
10\,T^*_1T_1 \,+\, 3\,T^*_2T_2 \,-\,T^*_3T_3 
\\[1ex]
on \quad \Lambda^3_{48}:&&
q(R) \,\,=\,\,
\tfrac{11}{4}\,T^*_1T_1 \,+\, \tfrac{15}{4}\,T^*_2T_2 \,+\, 
\tfrac{23}{4}\,T^*_3T_3 \,+\, \tfrac{7}{4}\,T^*_4T_4 \,-\,
\tfrac{1}{4}\,T^*_5T_5 \,-\, \tfrac{5}{4}T^*_6T_6
\\[1ex]
on \quad \Lambda^4_{27}:&&
q(R) \,\,=\,\,
\tfrac{7}{2}\,T^*_1T_1 \,-\, 2\,T^*_2T_2 
\\[1ex]
on \quad \Lambda^4_{35}:&&
q(R) \,\,=\,\,
6\,T^*_1T_1 \,+\, \tfrac{5}{2}\,T^*_2T_2 \,-\, 
 \tfrac{3}{2}\,T^*_4T_4 
\eea
\end{epr}

\section{Proof of the Theorems}

In this section we will prove Theorems~\ref{main1} and
\ref{main2} using the Weitzenb{\"o}ck formulas of 
Proposition~\ref{final1} and \ref{final2}. We will
first show that on manifolds with holonomy $\,\G_2\,$
and $\,\Spin_7\,$ any Killing form can be decomposed
into a sum of Killing forms belonging to the parallel
subbundles of the form bundle. Hence we may
assume that the Killing form is a section of
one of the irreducible components. Moreover we
only have to consider summands where the
endomorphism $q(R)$ acts non-trivially. Otherwise
the twistor form has to be parallel, as
follows from Corollary~\ref{ricci}.

It is easy to show that on Einstein manifolds the
codifferential of any twistor 2-form is either
zero or dual to a Killing vector field. However,
on a compact Ricci-flat manifold any Killing
vector field has to be parallel and it follows
after integration
that the twistor 2-form has to be coclosed.
This argument then implies the more general 
statement on twistor forms, contained in  
Theorems~\ref{main1} and
\ref{main2}.

\subsection{The holonomy decomposition}

Let $(M^n,g)$ be a manifold with holonomy $G$, which is assumed
to be a proper subgroup of $\SO_n$. In this situation the
bundle of p-forms decomposes into a sum of parallel subbundles,
$\Lambda^pTM = \oplus V_i$ and correspondingly, any p-form $u$
has a holonomy decomposition $u = \sum u_i$. If $u$ is a
twistor form, or even a Killing form, it remains in general 
not true for the components $u_i$. Nevertheless we have
such a property in the case of the exceptional holonomies.

\begin{elem}\label{holdeco}
Let $(M, g)$ be a compact manifold with holonomy $\G_2$ or
$\Spin_7$ and let $u$ be any form with holonomy decomposition
$u = \sum u_i$. Then $u$ is a Killing form or a $\ast$-Killing
form if and only if the same is true for all components $u_i$.
\end{elem}
\proof
We will use the characterization of Killing forms given in 
Proposition~\ref{integrabl}. Since the decomposition
$\Lambda^pTM = \oplus V_i$ is parallel, it is preserved by 
 $\nabla^*\nabla$ and $q(R)$. Thus for a Killing form $u$ all
the components satisfy the equation 
$\nabla^*\nabla u_i = \frac{1}{p}q(R)u_i$ and it remains to
verify whether the components are coclosed. 

Since manifolds with  holonomy $\G_2$ or $\Spin_7$ are
Ricci-flat we do not have to consider twistor 1-forms, which
are automatically parallel due to Corollary~\ref{ricci}. 
We start with the $\G_2$-case,
where it is enough to consider Killing forms.
The proof for $\ast$-Killing follows then from the duality
under the Hodge star operator.
Let $u = u_7 + u_{14}$ be the decomposition of a Killing 2-form
according to~(\ref{deco4}), then Corollary~\ref{ricci} shows that
$u_7$ is parallel and thus $u_{14} = u - u_7$ is coclosed.
In the case of a Killing 3-form we have the decomposition
$u = u_1 + u_7 + u_{27}$ and, as for 2-forms, it follows that
$u_1$ and $u_7$ have to be parallel, implying that $u_{27}$
is coclosed. The argument is the same for Killing forms
in degree 4 and 5, since the same representations are involved.

We now turn to the case of holonomy $\,\Spin_7$. 
In the case of twistor 4-forms it follows from the remark after
Proposition~\ref{integrabl} that any component is again
a twistor form. The further analysis will be given below.
The argument for forms in all other degrees is the
same as in the $\G_2$-case.

\subsection{Twistor forms on $G_2$-manifolds}

For $G_2$-manifolds we only have to consider two cases:
twistor forms in $\Lambda^2_{14}$ and in $\Lambda^3_{27}$.
In the first case we have three twistor operators (corresponding
to the decomposition~(\ref{deco})) and $T_3$ vanishes on
twistor forms, since the third summand in the 
decomposition~(\ref{deco}) of $T \otimes \Lambda^2_{14}$
belongs neither to the $\Lambda^1$ nor to $\Lambda^3$. 
In the second case we have five twistor operators
and  here $T_4$ and $T_5$ vanish on twistor forms,
since the last two summands in the
decomposition~(\ref{deco}) of $T \otimes \Lambda^3_{27}$
do not appear in the form spaces. From 
Proposition~\ref{final1} we have one Weitzenb{\"o}ck formula
for each case. However, to prove Theorem~\ref{main1}
we still need to compare the twistor operators with the
differential and codifferential.

Consider the case $\Lambda^2_{14}$, here the differential
splits as $d=d_7 + d_{27}$,~e.g. 
$d_7 = \sum (e_i \wedge \nabla_{e_i})_7$, with subscripts
denoting the projection onto the corresponding summand.
There is no part $d_1$, since the trivial representation
does not occur in the decomposition of $T \otimes\Lambda^2_{14} $.
The projection $\pr_1$ defining $T_1$ can be written in two ways:
$$
\pr_1: T \otimes  \Lambda^2_{14}
\stackrel{\pi_1}\longrightarrow
T
\stackrel{j_1}\longrightarrow
T \otimes  \Lambda^2_{14},\qquad
\pr_1: T \otimes  \Lambda^2_{14}
\stackrel{\pi_2}\longrightarrow
\Lambda^3_{7}
\stackrel{j_2}\longrightarrow
T \otimes  \Lambda^2_{14},
$$
with $\pi_1(X\otimes \alpha) = X \ins \alpha$ and
$\pi_2(X\otimes \alpha) = (X\wedge \alpha)_7$, and
the right inverses $j_1, j_2$. Hence, $du=0$
or $d^*u=0$ both imply $T_1u=0$. Similarly, 
$du=0$ implies $T_2u=0$.

Let $u$ be a $\ast$-Killing form in  $\Lambda^2_{14}$,
then $du=0$ implies $T_1u=0$ and $T_2u=0$. Hence all
twistor operators vanish on $u$ and the form has to
be parallel.
Let $u$ be a Killing form in  $\Lambda^2_{14}$, then
only the component $T_2u$ could be different from 0.
But the Weitzenb{\"o}ck formula of Proposition~\ref{final1}
and the equation $2\nabla^*\nabla u =q(R)u$ imply
$$
0 \,\,=\,\, 2\nabla^*\nabla u \,-\, q(R)u
\,\,=\,\ (2 - \tfrac{4}{3})T^*_2T_2u \ .
$$
Hence $T^*_2T_2u=0$, and after integration also $T_2u=0$,~i.e.
the form $u$ has to be parallel.

Consider now the case $\Lambda^3_{27}$, here the differential
splits as $d=d_7 + d_{27}$ and the codifferential as
$d^* = d^*_7 +  d^*_{14}$, again there is no component
$d_1$. With the same arguments as for $\Lambda^2_{14}$
we see that $du=0$ implies $T_1u=0$ and $T_2u=0$ and
$d^*u=0$ implies $T_1u=0$ and $T_3u=0$.
Indeed, the first to summands of the decomposition~(\ref{deco})
of $T\otimes \Lambda^3_{27}$ are also summands
of the 4-forms. Similarly, the first and the third
summand are components of the 2-forms.

Let $u$ be Killing form in $\Lambda^3_{27}$, then only
$T_2u$ could be different from zero. But the
Weitzenb{\"o}ck formula and the equation 
$3\nabla^*\nabla u =q(R)u$ implies $(3-2)T^*_2T_2u=0$ and
$u$ again has to be parallel.
Finally, in the case of a $\ast$-Killing form $u$ in 
$\Lambda^3_{27}$, we have to use the Weitzenb{\"o}ck 
formula and the equation  $4\nabla^*\nabla u =q(R)u$
to show the vanishing of $T_3u$.
\qed

\begin{ere}
Using Lemma~\ref{adjoint} and explicit expressions for
the projections onto the irreducible components of
the form bundle it is possible to determine
the precise relation between the operators
$T^*_iT_i$ and similar operators in terms of
the components of $d$ and $d^*$,~e.g. on 
the bundle $\Lambda^2_{14}$ one finds:
$\,d^*_7 d_7 = T^*_1T_1,\, d d^* = 4T^*_1T_1\,$
and $\,3d^*_{27}d_{27}=7T^*_2T_2$.
\end{ere}

\subsection{Twistor forms on $\Spin_7$-manifolds}

Since $q(R)$ acts trivially on the representations
$\C,\,\Lambda^2_7\,$ and $\,T$, 
any twistor form in one of these components is automatically
parallel. Hence, we only have to consider twistor forms in
the subbundles
$\,\Lambda^2_{21},\, 
\Lambda^3_{48}$,\,
$\,\Lambda^4_{27}$\, and $\,\Lambda^4_{35}$.
The argument is now similar to the $\G_2$-case. For any
Killing or $\ast$-Killing form $u$ in one of these bundles
we show that all twistor operators vanish on $u$, such
that the form has to be parallel. 
%We use the numbering
%of the twistor operator according to the decomposition~(\ref{deco6}).

Let $u \in\Lambda^2_{21} $ be a twistor 2-form, then $T_3u=0$,
since the third summand in the decomposition~(\ref{deco6})
of $T\otimes \Lambda^2_{21}$ belongs neither to $\Lambda^1$
nor to $\Lambda^3$. The representation $T \cong \Lambda^1$ appears 
also as summand in the 3-form. 
Hence, as in the $\G_2$-case, we 
see that $du=0$ implies $T_1u = T_2u=0$ and $d^*u=0$
implies $T_1u=0$. Thus $\ast$-Killing forms are automatically
parallel. Let $u$ be a Killing 2-form, then
$2\nabla^*\nabla u= q(R)u$ and the Weitzenb{\"o}ck formula
of Proposition~\ref{final2}
imply $\,(2-3)T^*_2T_2u\,$ and $\,T_2u=0\,$. Thus also
on Killing forms all twistor operators vanish.

Let $u \in\Lambda^3_{48} $ be a twistor 3-form,~i.e. $T_5u=T_6u=0$.
The representation $\Lambda^2_{7}$,~i.e. the third summand
in the decomposition~(\ref{deco6}) of $T\otimes \Lambda^3_{48}$,
appears as summand in the 2- and 4-forms.
Hence,
$du=0$ implies $T_1u =T_3u=T_4u=0$ and $d^*u=0$ implies
$T_2u=T_3u=0$. Let $u$ be a $\ast$-Killing form, then
$5\nabla^*\nabla u= q(R)u$ and we find 
$(5 - \frac{15}{4})T^*_2T_2u=0$, thus $T_2u=0$. Let $u$ be
a Killing form, then $3\nabla^*\nabla u= q(R)u$ and
$\,(3-\frac{11}{4})T^*_1T_1 +(3-\frac{7}{4})T^*_4T_4u=0\,$
implies $T_1u=T_4u=0$.

Let $u\in \Lambda^4_{27}$ be a twistor 4-form,~i.e. $T_2u=0$.
Then $4\nabla^*\nabla u= q(R)u$ implies 
$\,(4-\frac{7}{2})T^*_1T_1u=0\,$ and thus $\,T_1u=0$.
Hence, any twistor 4-form in $\Lambda^4_{27}$ has to be
parallel, which completes the proof Lemma~\ref{holdeco}.
We had already seen that any component of a twistor 
4-form is again a twistor form. Now we see that three of
the four components have to be parallel. Hence, for a
Killing form all components are again coclosed,~i.e.
again Killing forms.

Finally, let $u\in \Lambda^4_{35}$ be a twistor 4-form,~i.e.
$T_3u=T_4u=0$. Since both remaining summands
$\,T\, $ and $\,V^b_{112}\,$ are subbundles
of $\Lambda^3\cong \Lambda^5$ we see that $T_1$ and $T_2$
vanish for closed or coclosed twistor forms, thus they
have to be parallel.
\qed

\begin{ere}
For the representation $\,\Lambda^4_{35}\,$ we find  in 
\cite{gu} an additional Weitzenb{\"o}ck 
formula (with vanishing curvature term). This equation
then shows that indeed any twistor 4-form on a 
$\Spin_7$-manifold has to be parallel.
\end{ere}

%%%%%%%%%%%%%%%%%%%%%%%%%%%%%%%%%%%%%%%%%%%%%%%%%%%%%%%%%%%%%%%%%%%%%%%%%

\labelsep .5cm


\begin{thebibliography}{22}
{\footnotesize

\bibitem{pg3}
  \textsc{Apostolov,~V., Calderbank,~D., Gauduchon,~P.},
  \textit{Hamiltonian $2$--forms in K{\"a}hler geometry~I},
  \textrm{math.DG/0202280 (2002)}.


\bibitem{fau} 
\textsc{Belgun,~F., Moroianu,~A.,  Semmelmann,~U.},
{\sl Killing forms on symmetric spaces,}
{\textrm math.DG /0409104, (2004)}.


 \bibitem{be1}
  \textsc{Benn, I. M., Charlton, P., Kress, J. },
  \textit{Debye potentials for Maxwell and Dirac fields from a generalization 
    of the   Killing-Yano equation},
  \textrm{J. Math. Phys. 38 (1997), no. 9, 4504--4527}.


 \bibitem{be2}
  \textsc{Benn, I. M., Charlton, P.},
  \textit{Dirac symmetry operators from conformal Killing-Yano tensors},
  \textrm{ Classical Quantum Gravity 14 (1997), no. 5, 1037--1042}.


\bibitem{pg2} \textsc{Calderbank, D., Gauduchon, P., Herzlich, M.},
   {\sl Refined Kato Inequalities and Conformal weights in
   Riemannian Geometry,}
   \textrm{J. Funct. Anal. {\bf 173} (2000), 214--255. }


 \bibitem{pg1} \textsc{Gauduchon, P.,}
  {\sl  Structures de Weyl et 
theoremes d'annulation sur une variete conforme autoduale, }
 \textrm{Ann. Scuola Norm. Sup. Pisa Cl. Sci. (4)  {\bf 18}  (1991),  no. 4, 563--629.}

\bibitem{homma}
  \textsc{Homma,~Y.},
  \textit{
Bochner-Weitzenb{\"o}ck formulas and curvature actions on Riemannian manifolds},
  \textrm{math.DG/0307022 }.


\bibitem{au1} \textsc{Moroianu, A., Semmelmann, U.},
{\sl Twistor forms on K{\"a}hler manifolds, }
\textrm{Ann. Scuola Norm. Sup. Pisa Cl. Sci. (4) {\bf II} (2003), 823--845 }


\bibitem{au2} \textsc{Moroianu,~A., Semmelmann,~U.}
{\sl Killing forms on quaternion K{\"a}hler manifolds,}
{\textrm Math.DG /0403242 (2004)}


\bibitem{prod} \textsc{Moroianu,~A., Semmelmann,~U.} {\sl Twistor
 Forms on Riemannian Products, } \textrm{math.DG/0407063 (2004)}.




 \bibitem{penrose}
  \textsc{Penrose, R., Walker, M.},
  \textit{On quadratic first integrals of the geodesic equations for 
  type $\{22\}$ spacetimes,  }
  \textrm{Comm. Math. Phys. 18 1970 265--274.  }


 \bibitem{uwe} \textsc{Semmelmann,~U.,}
  {\sl Conformal Killing forms on Riemannian manifolds, }
  \textrm{Math. Z. {\bf 243} (2003), 503--527.}


 \bibitem{gu} \textsc{Semmelmann, U., Weingart, G.},
  {\sl Weitzenb{\"o}ck formulas for manifolds with special holonomy, }
  in preparation.

 \bibitem{yama3}
  \textsc{Yamaguchi, S.},
  {\sl On a Killing $p$-form in a compact K{\"a}hlerian manifold},
  Tensor (N.S.) {\bf 29} (1975), no. 3, 274--276. 



 \bibitem{yano}
  \textsc{Yano, K.},
  \textit{Some remarks on tensor fields and curvature, }
  \textrm{Ann. of Math. (2) 55, (1952). 328--347.}



}
\end{thebibliography}
\end{document}